\documentclass[11pt]{article}
\usepackage{amsmath}
\usepackage{amsthm}
\usepackage{epsfig}
\usepackage{leftidx}
\usepackage{graphicx,subfigure}
\usepackage{amsfonts}
\usepackage{graphicx, color, epstopdf}
\usepackage{cite}
\usepackage{amssymb,amsmath}
\usepackage{multicol}
\usepackage{booktabs}
\usepackage[center]{caption2}
\usepackage{multirow}
\usepackage{enumerate}
\usepackage{enumitem}
\usepackage{algpseudocode,algorithm,algorithmicx}

\newtheorem{theorem}{Theorem}[section]
\newtheorem{lemma}{Lemma}[section]

\newcommand{\defeq}{:=}
\newcommand{\zd}{\,\mathrm{d}}
\newcommand{\abs}[1]{\left|#1\right|}
\newcommand{\absb}[1]{\big|#1\big|}

\newcommand{\bra}[1]{\left(#1\right)}
\newcommand{\brab}[1]{\big(#1\big)}
\newcommand{\braB}[1]{\Big(#1\Big)}

\newcommand{\kbra}[1]{\left[#1\right]}
\newcommand{\kbrab}[1]{\big[#1\big]}

\newcommand{\myinnerb}[1]{\big\langle#1\big\rangle}

\newcommand{\mynormb}[1]{\big\|#1\big\|}

\newcommand{\mynormt}[1]{\|#1\|}

\title{Discrete energy analysis of the
	 third-order variable-step BDF time-stepping for diffusion equations}
\author{Hong-lin Liao\thanks{ORCID 0000-0003-0777-6832. School of Mathematics,
	Nanjing University of Aeronautics and Astronautics,
Nanjing 211106, China; Key Laboratory of Mathematical Modelling
and High Performance Computing of Air Vehicles (NUAA), MIIT, Nanjing 211106, China.
Emails: liaohl@nuaa.edu.cn and liaohl@csrc.ac.cn.
This author's work is supported by NSF of China
under grant number 12071216.}
\quad Tao Tang\thanks{Division of Science and Technology,
    BNU-HKBU United International College, Zhuhai, Guangdong Province,China.
    Email: tangt@sustech.edu.cn.
    This author's work is partially supported by NSF of China
    under grant numbers 11731006 and NNW2018-ZT4A06 project.}
\quad Tao Zhou\thanks{Institute of Computational Mathematics
    and Scientific/Engineering Computing,
    Academy of Mathematics and Systems Science,
    Chinese Academy of Sciences, Beijing, 100190, China.
    Email: tzhou@lsec.cc.ac.cn.
    This author's work is partially supported by NSF of China
    (under grant numbers 11822111, 11688101, and 11731006).}
}
\date{}

\begin{document}

\maketitle

\begin{abstract}
This is one of our series works on discrete energy analysis of the
variable-step BDF schemes. In this part, we present stability and convergence analysis of the
third-order BDF (BDF3) schemes with variable steps for linear diffusion equations,
see e.g. [SIAM J. Numer. Anal., 58:2294-2314] and [Math. Comp., 90: 1207-1226] for
our previous works on the BDF2 scheme. To this aim,
we first build up a discrete gradient structure of the variable-step BDF3 formula under the condition that
the adjacent step ratios are less than 1.4877, by which we can establish a discrete energy dissipation law.
Mesh-robust stability and convergence analysis in the $L^2$ norm are then obtained. Here the mesh robustness means that
the solution errors are well controlled by the maximum time-step size but
independent of the adjacent time-step ratios. We also present numerical tests to support our theoretical results.
\\[1ex]
\\[1ex]
\end{abstract}

\section{Introduction}
\setcounter{equation}{0}

In this paper, we aim to develop a discrete energy technique
for the stability and convergence of three-step
backward differentiation formula (BDF3) with variable time-steps.
To this end, we consider the linear reaction-diffusion problem in a bounded convex domain $\Omega$,
\begin{align}\label{eq: diffusion problem}
	\partial_tu-\varepsilon\Delta u&=\kappa(x)u+f(t,x)\quad
	\text{for $x\in\Omega$ and $0<t<T$,}
\end{align}
subject to the Dirichlet boundary condition~$u=0$ on  a smooth boundary~$\partial\Omega$, and the initial data $u(0,x)=u_0$ for $x\in\Omega$.
We assume that the diffusive coefficient $\varepsilon>0$ is a constant
and the reaction coefficient $\kappa(x)$ is smooth  and bounded by $\kappa^{*}>0$.

The BDF schemes are widely used
for stiff or differential-algebraic problems
\cite{HairerNorsettWanner:1992,HairerWanner:2002}.
Recently, they were also applied for simulating hyperbolic systems with multiscale
relaxation \cite{AlbiDimarcoPareschi:2020} and stiff kinetic equations \cite{DimarcoPareschi:2017}.
For such applications, BDF schemes with variable steps are shown to be computationally efficient in capturing the multi-scale time behaviors
\cite{CalvoGrandeGrigorieff:1990,ChenWangYanZhang:2019,DeCariaGuzelLaytonLi:2021,
	HairerWanner:2002,LiaoJiZhang:2020pfc,
	AlbiDimarcoPareschi:2020,WangRuuth:2008}.
However, rigorous theoretically analysis (stability and convergence) for variable-step
BDF schemes is challenging. This motivates our serious works on this topic, and one can find our previous works on the variable-step BDF2 scheme \cite{LiaoZhang:2021,LiaoTangZhou:2020bdf2}.

To begin, we consider the temporal mesh
$$0=t_0<t_1<t_2<\cdots<t_N=T$$
with the variable time-step
$$\tau_k\defeq t_k-t_{k-1}, \quad 1\le k\le N.$$
The maximum step size and the adjacent time-step ratios are defined respectively as
$$\tau\defeq\max_{1\le k\le N}\tau_k, \quad  r_{k}\defeq\frac{\tau_{k}}{\tau_{k-1}} \,\,\, \textmd{for}\,\,\,  2\le k\le N.$$
For any sequences $\{v^n\}_{n=0}^N$, we denote $\partial_{\tau}v^n:=(v^n-v^{n-1})/\tau_n$.
Then by taking $v^n=v(t_n)$, the variable-step BDF3 formula \cite{CalvoGrigorieff:2002} yields
\begin{align}\label{eq: BDF3 formula}
	D_3v^n=&\,d_0(r_n,r_{n-1})\partial_{\tau} v^{n}
	+d_1(r_n,r_{n-1})\partial_{\tau} v^{n-1}+d_2(r_n,r_{n-1})\partial_{\tau} v^{n-2},
\end{align}
where
\begin{align}
d_0(x,y):=&\,\frac{1+2x}{1+x}+\frac{xy}{1+y+xy},\label{funs: d0}\\
d_1(x,y):=&\,-\frac{x}{1+x}
-\frac{xy}{1+y+xy}
-\frac{xy^2}{1+y+xy}\frac{1+x}{1+y},\label{funs: d1}\\
d_2(x,y):=&\,\frac{xy^2}{1+y+xy}
\frac{1+x}{1+y},\qquad\text{for $x, y\ge0$.}\label{funs: d2}
\end{align}
Without losing the generality, we assume that the discrete solution $u^1$ and $u^2$ are given. We now consider a time-discrete solution for the diffusion equations, $u^k(x)\approx u(t_k,x)$ for~$x\in\Omega$, by the following variable-step BDF3 time-stepping scheme
\begin{align}\label{eq: time-discrete IBVP}
	D_3u^{k}=\varepsilon\Delta u^{k}
	+\kappa u^{k}+f^k,\quad, 3\le k\le N,
\end{align}
where $f^k(x)=f(t_k,x)$.

To the best of our knowledge, there are very few theoretical results
on the variable-step BDF3 scheme in literature.
For linear diffusion problems,
Calvo and Grigorieff \cite{CalvoGrigorieff:2002} established the following $L^2$ norm stability estimate
under the step-ratio condition $r_k<1.199$,
\begin{align*}
	\mynormb{u^{n}}
	\le C\exp\bra{C\Gamma_n}\braB{\mynormb{u_{0}}+\sum_{j=1}^{n}\tau_j\mynormb{f^j}}\quad\text{for $n\ge1$,}
\end{align*}
where the ratio-dependent \textit{prefactor} $\Gamma_n$ is given by
$$\Gamma_n:=\sum_{k=2}^{n}\absb{r_{k}-r_{k-1}}.$$
This means that any mixture of the $k$-step BDF $(k\in\{1,2,3\})$ is stable provided
the number of changes between increasing and decreasing mesh-sizes is uniformly bounded.
Nonetheless, the quantity $\exp\bra{C\Gamma_n}$ can be unbounded
since $\Gamma_n$ could blow up for certain time-step
series at vanishing time-steps. To see this, we consider a specific class of time-steps
$\{\cdots, \tau_1,\mu\tau_1,\tau_1,\mu\tau_1,\cdots\}$ with a positive constant $\mu\neq1$.
Then for a finite time $T=\frac{M}{2}(1+\mu)\tau_1$, one has
$$\Gamma_M=\sum_{k=2}^{M}\absb{\mu-\mu^{-1}}
=(M-1)\absb{\mu-\mu^{-1}}\rightarrow\infty\qquad\text{as $\tau_1\rightarrow0$}.$$

To remedy this issue, we present shall in this work a new framework for analyzing the variable step BDF3 scheme. Our contribution is two folds:
\begin{itemize}
  \item We build up a discrete graident structure of BDF3 formula
  under the step-ratio condition $0<r_k<R_{e}$,
  where $R_{e}\approx1.4877$
  is the unique positive root of
  $$d_1(R_{e},0)+\tfrac{7}{10}\sqrt{R_{e}}d_2(R_{e},R_{e})=0.$$
  Based on this, we present the discrete energy stability  for the variable-step BDF3 scheme.

  \item We further present the stability
  and convergence analysis in the $L^2$ norm for the BDF3 scheme under the adjacent step ratios $0<r_k<R_{e}.$ We also show that these analysis results are mesh-robustly, which means that the associated \textit{prefactors} in our analysis are independent of the time-step ratios $r_k$.
  That is, unbounded quantities such as $\Gamma_n$ in \cite{CalvoGrigorieff:2002} are removed.
\end{itemize}
The rest of this work is organized as following. In the next section, we provide with some preliminary tools, e.g., discrete orthogonal convolution kernels, for our analysis. In section 3, we present the energy stability of the variable-step BDF3 scheme. The stability and convergence analysis of the variable-step BDF3 scheme is then presented in section 4. This is followed by some numerical examples in section 5.

\section{Discrete orthogonal convolution kernels}
In this section, we present some preliminary tools for our analysis.
We remark that the mesh-robust stability and convergence of the variable-step
BDF2 scheme have been established in our previous works \cite{LiaoJiZhang:2020pfc,LiaoZhang:2021}.
Nonetheless, the extension to BDF3 formula is theoretically challenging
due to the additional degrees of freedom. This work is builds upon our previous work \cite{LiLiao:2021},
where the stability of variable-step BDF3 method for nonlinear ODE problems
was verified under the step-ratio condition $r_k<R_3\approx2.553$.

To begin, we write the BDF3 formula \eqref{eq: BDF3 formula} as
a convolution of local difference quotients
\begin{align}\label{eq: alter BDFk formula}
	D_{3}v^n:=\sum_{j=1}^nd^{(n)}_{n-j}\partial_{\tau}  v^j\quad\text{for $n\ge3$,}
\end{align}
where the associated discrete BDF3 kernels
\begin{align}\label{eq: BDF3 kernels}
	&d^{(n)}_j:=d_j(r_n,r_{n-1})\;\;\text{for $j=0,1,2$,}\quad\text{and}\quad	d^{(n)}_j:=0
	\quad\text{for $n\ge j+1\ge 4$}.
\end{align}
Assume always that the summation $\sum_{k=i}^{j}\cdot$ to be zero
and the product $\prod_{k=i}^{j}\cdot$ to be one if the index $i>j$.
As for the BDF3 kernels $d^{(n)}_{n-j}$ with any fixed indexes $n$,
we recall a class of discrete orthogonal convolution
(DOC) kernels $\big\{\vartheta_{n-j}^{(n)}\big\}_{j=3}^n$ by a recursive procedure,
also see \cite{LiaoZhang:2021},
\begin{align}\label{eq: BDFk-DOC procedure}
	\vartheta_{0}^{(n)}:=\frac{1}{d^{(n)}_{0}}\quad\text{and}\quad
	\vartheta_{n-j}^{(n)}:=-\frac{1}{d^{(j)}_{0}}
	\sum_{i=j+1}^{n}\vartheta_{n-i}^{(n)}d^{(i)}_{i-j}\quad\text{for $3\leq j\le n-1$.}
\end{align}
Obviously, the DOC kernels $\vartheta_{n-j}^{(n)}$ satisfy the following discrete orthogonality identity
\begin{align}\label{eq: BDFk-DOC identity}
	\sum_{i=j}^{n}\vartheta_{n-i}^{(n)}d^{(i)}_{i-j}\equiv\delta_{nj}
	\quad\text{for any $3\leq j\le n$,}
\end{align}
where $\delta_{nj}$ is the Kronecker delta symbol with $\delta_{nj}=0$ if $j\neq n$.
Furthermore, with the
identity matrix $I_{m\times m}$ ($m:=n-2$), the
above discrete orthogonality identity \eqref{eq: BDFk-DOC identity} also implies
$$\Theta_3D_3=I_{m\times m},$$
where the two $m\times m$ matrices $D_{3}$
and $\Theta_3$ are defined by
\[
D_{3}:=
\left(
\begin{array}{cccc}
	d_{0}^{(3)}  &                  &  & \\
	d_{1}^{(4)}  &d_{0}^{(4)}  &  & \\
	\vdots           &\vdots           &\ddots  &\\
	d_{n-3}^{(n)}&d_{n-2}^{(n)}&\cdots  &d_{0}^{(n)}  \\
\end{array}
\right)\;\;\text{and}\;\;
\Theta_3:=
\left(
\begin{array}{cccc}
	\vartheta_{0}^{(3)}  &                  &  & \\
	\vartheta_{1}^{(4)}  &\vartheta_{0}^{(4)}  &  & \\
	\vdots           &\vdots           &\ddots  &\\
	\vartheta_{n-3}^{(n)}&\vartheta_{n-2}^{(n)}&\cdots  &\vartheta_{0}^{(n)}  \\
\end{array}
\right).
\]
Obviously, one has
$D_{3}\Theta_3=I_{m\times m},$
which implies the following mutual orthogonality identity
\begin{align}\label{eq: mutual BDFk-DOC identity}
	\sum_{i=j}^{n}d_{n-i}^{(n)}\vartheta^{(i)}_{i-j}\equiv\delta_{nj}
	\quad\text{for any $3\leq j\le n$.}
\end{align}

Note that,
this identity \eqref{eq: mutual BDFk-DOC identity} will be used to
study the discrete property of DOC kernels; while
the above identity \eqref{eq: BDFk-DOC identity} will be used to reformulate
the discrete scheme \eqref{eq: time-discrete IBVP}.
By exchanging the summation order and using \eqref{eq: BDFk-DOC identity}, one has
\begin{align}\label{eq: BDFk-DOC transform}
	\sum_{i=3}^{n}\vartheta_{n-i}^{(n)}D_{3}v^i
	=&\,\sum_{i=3}^{n}\vartheta_{n-i}^{(n)}
	\sum_{j=1}^{2}d^{(i)}_{i-j}\partial_{\tau}  v^j
	+\sum_{i=3}^{n}\vartheta_{n-i}^{(n)}
	\sum_{j=3}^id^{(i)}_{i-j}\partial_{\tau}  v^j\nonumber\\
	=&\,\sum_{j=1}^{2}\partial_{\tau}v^j
	\sum_{i=3}^n\vartheta_{n-i}^{(n)}d^{(i)}_{i-j}+
	\sum_{j=3}^{n}\partial_{\tau}v^j
	\sum_{i=j}^n\vartheta_{n-i}^{(n)}d^{(i)}_{i-j}\nonumber\\
	=&\,\mathcal{I}_{3}^n[v]+\partial_{\tau}  v^n\qquad\text{for $n\ge3$,}
\end{align}
where $\mathcal{I}_{3}^n[v]$ represents the starting effect on the numerical solution at $t_n$, or
\begin{align}\label{eq: initial effect BDF3-DOC}
	\mathcal{I}_{3}^n[v]:=&\,\sum_{j=1}^{2}\partial_{\tau}v^j
	\sum_{i=3}^n\vartheta_{n-i}^{(n)}d^{(i)}_{i-j}
	=\partial_{\tau}v^2
	\sum_{i=3}^n\vartheta_{n-i}^{(n)}d^{(i)}_{i-2}
	+\vartheta_{n-3}^{(n)}d^{(3)}_{2}\partial_{\tau}v^1
	\quad\text{for $n\ge3$.}
\end{align}
Multiplying both sides of the equation \eqref{eq: time-discrete IBVP} by the DOC kernels
$\vartheta_{n-k}^{(n)}$, and summing $k$ from $3$ to $n$,
we apply \eqref{eq: BDFk-DOC transform} to get the following equivalent form
\begin{align}\label{eq: DOC-BDF3 equation}
	\partial_{\tau}u^n
	=-\mathcal{I}_3^n[u]+\sum_{k=3}^n\vartheta_{n-k}^{(n)}\bra{\varepsilon\Delta u^k+\kappa u^k}
	+\sum_{k=3}^n\vartheta_{n-k}^{(n)}f^k
	\quad\text{for $3\le n\le N$.}
\end{align}
This formulation expresses the BDF3 solution at time $t_n$ as
a (global) convolution summation of all previous solutions with DOC kernels
$\vartheta_{n-k}^{(n)}$ as the convolutional weights.

Next section constructs a discrete gradient structure of variable-step BDF3 formula
and derives the energy stability of BDF3 scheme
via the original form \eqref{eq: time-discrete IBVP}.
Section 3 addresses the $L^2$ norm stability and convergence analysis
via the discrete convolution form \eqref{eq: DOC-BDF3 equation}.
Some numerical examples are presented
in the last section to support our theoretical results.

\section{Positive definiteness and energy stability}
\setcounter{equation}{0}

To show the energy stability of the variable-step BDF3 scheme, we first investigate sufficient conditions on the adjacent time-step ratios $r_k$ so that the discrete kernels $\big\{\tau_nd_{n-k}^{(n)}\big\}$ are positive definite.

For certain adaptive time-stepping process, one may choose the step size $\tau_{m+1}$
(or the step ratio $r_{m+1}$) properly according to the information
from previous time-step ratios $\{r_k\}_{k=2}^{m}$.
Actually, the positive definiteness should be determined by the
eigenvalues of the pentadiagonal symmetric matrix $B_3=B_L+B_L^T$, where
\begin{align}\label{lemproof: BDF3 lower triangle matrix}
	B_L:=
	\left(
	\begin{array}{ccccc}
		\tau_3d_0^{(3)} & && &\vspace{0.15cm}\\
		\tau_4d_1^{(4)}&\tau_4d_0^{(4)} &&&\vspace{0.15cm}\\
		\tau_5d_2^{(5)}&\tau_5d_1^{(5)} &\tau_5d_0^{(5)}&&\vspace{0.15cm}\\
		&\ddots  &\ddots&\ddots&\vspace{0.1cm}\\
		&&\tau_nd_2^{(n)}&\tau_nd_1^{(n)}&\tau_nd_0^{(n)}\\
	\end{array}
	\right).
\end{align}
A sufficient and necessary condition for the positive definiteness of $B_3$ would be a
certain combination involving all time-step ratios; however,
it remains open at this moment. We consider only certain
restriction of each step ratio, like $0<r_k<R_{e}$
for a fixed positive constant $R_{e}<R_3\approx2.553$,
the recent stability constraint \cite{LiLiao:2021} for the ODE problems.

\begin{table}[htb!]
	\begin{center}
		\caption{Minimum eigenvalue of step-rescaled matrix
			$\widetilde{B}_3$ on random meshes.}
		\label{exam: Minimum eigenvalue}\vspace*{0.3pt}
		\def\temptablewidth{0.7\textwidth}
		{\rule{\temptablewidth}{0.5pt}}
		\begin{tabular*}{\temptablewidth}{@{\extracolsep{\fill}}ccccc}
			$n$    &$R_e=1.20$	&$R_e=1.50$ &$R_e=1.69$	&$R_e=1.70$	\\
			\hline
			50        &1.12 &5.08e-01	&6.12e-02  &-4.55e-02     	\\
			100       &1.07	&4.35e-01	&4.58e-02  &-5.29e-02		\\
			200       &1.08	&4.18e-01	&-2.06e-02 &-8.49e-02		\\
		\end{tabular*}
		{\rule{\temptablewidth}{0.5pt}}
	\end{center}
\end{table}	


Numerical tests on random meshes
are performed to examine the positive definiteness of
$B_3$  via the step-rescaled matrix
$\widetilde{B}_3:=\Lambda_{\tau}^{-1}(B_L+B_L^T)\Lambda_{\tau}^{-1}$,
where $\Lambda_{\tau}=\mathrm{diag}(\sqrt{\tau_3},\sqrt{\tau_4},\cdots,\sqrt{\tau_n})$.
We take a finite time $T=1$ with $n$ grid points and
let $r_k$ $(2\le k\le n)$ be uniformly distributed random numbers over $(0,R_e)$.
Table \ref{exam: Minimum eigenvalue} lists
the minimum eigenvalue
(each data is the minimum value of 200 runs on different random meshes)
of  $\widetilde{B}_3$ for the fixed step-ratio limits
$R_e=1.20$, $1.50$, $1.69$ and $1.70$.
To ensure the positive definiteness, Table \ref{exam: Minimum eigenvalue} suggests that the
maximum step-ratio limit $R_e<1.69$ is necessary, while
we will prove theoretically that $R_e<1.4877$ is sufficient in the next subsection.

\subsection{Discrete gradient structure}
To derive the energy stability of numerical scheme, we need
a discrete gradient structure of variable-step BDF3 formula.
In the following, we will seek two nonnegative quadratic
functionals $G$ and $F$ such that
\begin{align}\label{eq: BDF3 energy quantity}
	J_n:=2v_n\tau_n\sum_{j=3}^nd_{n-j}^{(n)}v_j
	=G[v_n,v_{n-1}]-G[v_{n-1},v_{n-2}]+F[v_n,v_{n-1},v_{n-2}],
\end{align}
for $n\ge 3$, where the discrete BDF3 kernels $d_{j}^{(n)}$ are defined by \eqref{eq: BDF3 kernels}, such that
the associated quadratic form is positive definite
\begin{align*}
	2\sum_{k=3}^nv_k\tau_k\sum_{j=3}^kd_{k-j}^{(k)}v_j
	>0\quad\text{if $v_k\not\equiv0$.}
\end{align*}
This seems to be a difficult task due to the presence of variable kernels $d_{j}^{(n)}$,
refer to the recent comments in \cite[section 3.4]{DeCariaGuzelLaytonLi:2021}.
For the uniform case with $r_k\equiv1,$ it has been shown in \cite{Pierre:2021} that
the uniform BDF3 formula admits the following discrete gradient structure
\begin{align*}
	2v_n\braB{\frac{11}{6}v_n-\frac{7}{6}v_{n-1}+\frac{1}{3}v_{n-2}}
	=&\,\frac{1}{3}v_{n}^2+\frac1{3}(\tfrac7{4}v_n-v_{n-1})^2
	-\kbra{\frac{1}{3}v_{n}^2+\frac1{3}(\tfrac7{4}v_n-v_{n-1})^2}\\
	&\,+\frac{95}{48}v_n^2+\frac{1}{3}(v_n-\tfrac{7}{4}v_{n-1}+v_{n-2})^2.
\end{align*}
This decomposition is optimal in the sense that the minimum eigenvalue bound $\frac{95}{48}$
of the associated quadratic form is sharp, see  \cite[Lemma 2.4]{LiaoTangZhou:2021bdf345},
due to the Grenander-Szeg\"{o} theorem.

To deal with the variable-step case, our first task is to introduce a step-rescale transform $v_j=w_j/\sqrt{\tau_j}$ for $j\ge1$,
cf. the above step-rescaled matrix
$\widetilde{B}_3$, to remove the time-step factor $\tau_n$ in the discrete kernels
of \eqref{eq: BDF3 energy quantity}. One can get
\begin{align}\label{eq: modified BDF3 kernels}
	J_n=2w_n\sum_{j=3}^n
	\tilde{d}_{n-j}^{(n)}w_j\quad
	\text{with}\quad
	\tilde{d}_{n-j}^{(n)}:=\frac{\sqrt{\tau_n}}{\sqrt{\tau_j}}d_{n-j}^{(n)}
	\quad \text{for $n\ge j+2\ge 3$}.
\end{align}
It is reasonable to assume that
 \begin{align*}
 	\widetilde{G}[w_n,w_{n-1}]:=&\,a_{n+1}w_n^2+b_{n+1}(\gamma w_n-w_{n-1})^2,\\
 	\widetilde{F}[w_n,w_{n-1},w_{n-2}]:=&\,
 	p_{n+1}w_n^2+q_{n+1}(\gamma w_n-w_{n-1})^2+c_n\brab{w_n-\gamma w_{n-1}+w_{n-2}}^2,
 \end{align*}
where the nonnegative variable coefficients $a_n,b_n,c_n,p_n,q_n$ and the real parameter $\gamma$ (for which $\gamma=7/4$ would not necessarily be optimal) is determined such that
 \begin{align}\label{eq: BDF3 decomposition assumption}
	J_n=\widetilde{G}[w_n,w_{n-1}]
	-\widetilde{G}[w_{n-1},w_{n-2}]
	+\widetilde{F}[w_n,w_{n-1},w_{n-2}].
\end{align}
The principle of identity gives the following relationships for the undetermined coefficients:
\begin{align*}
	\text{coefficients of $w_{n}w_{n-2}$:}\quad&\, 2c_{n}=2\tilde{d}_{2}^{(n)};\\
	\text{coefficients of $w_{n-2}^2$ ($w_{n-1}w_{n-2}$):}\quad&\,
	-b_{n}+c_n=0;\\
	\text{coefficients of $w_{n}w_{n-1}$:}\quad&\,
	-2\gamma(b_{n+1}+q_{n+1})-2\gamma c_n=2\tilde{d}_{1}^{(n)};\\
	\text{coefficients of $w_{n-1}^2$:}\quad&\,
	-a_{n}+(b_{n+1}+q_{n+1})-\gamma^2b_{n}+\gamma^2c_n=0;\\
	\text{coefficients of $w_{n}^2$:}\quad&\,	a_{n+1}+p_{n+1}+\gamma^2(b_{n+1}+q_{n+1})+c_n
	=2\tilde{d}_{0}^{(n)}.
\end{align*}
They yield that
\begin{align*}
	& a_{n}=-\frac1{\gamma}
	\tilde{d}_{1}^{(n)}
	-\tilde{d}_{2}^{(n)},\quad
	b_{n}=\tilde{d}_{2}^{(n)},\quad
	c_{n}=\tilde{d}_{2}^{(n)},\\
	&p_{n+1}=2\tilde{d}_{0}^{(n)}
	+\gamma\tilde{d}_{1}^{(n)}
	+(\gamma^2-1)\tilde{d}_{2}^{(n)}
	+\frac1{\gamma}\tilde{d}_{1}^{(n+1)}
	+\tilde{d}_{2}^{(n+1)},\\
	&q_{n+1}=a_n-b_{n+1}=-\frac1{\gamma}\tilde{d}_{1}^{(n)}
	-\tilde{d}_{2}^{(n)}
	-\tilde{d}_{2}^{(n+1)}.
\end{align*}
According to the definitions in \eqref{eq: modified BDF3 kernels} and \eqref{eq: BDF3 kernels},
the coefficients $b_n$ and $c_n$ are always positive, while $a_n$ is also positive if $q_{n+1}\ge0$.
Thus the above assumption \eqref{eq: BDF3 decomposition assumption} requires
\begin{align}
	q_{n+1}=&\,-\frac1{\gamma}\tilde{d}_{1}^{(n)}
	-\tilde{d}_{2}^{(n)}
	-\tilde{d}_{2}^{(n+1)}\ge0,
	\label{eq: BDF3 decomposition condition1}\\
	p_{n+1}=&\,2\tilde{d}_{0}^{(n)}-\tilde{d}_{2}^{(n)}
	+\gamma^2\braB{\frac1{\gamma}\tilde{d}_{1}^{(n)}+\tilde{d}_{2}^{(n)}}	
	+\frac1{\gamma}\tilde{d}_{1}^{(n+1)}
	+\tilde{d}_{2}^{(n+1)}>0\quad \text{for all $r_k<R_{e}$.}
	\label{eq: BDF3 decomposition condition2}
\end{align}

The inequality system
\eqref{eq: BDF3 decomposition condition1}-\eqref{eq: BDF3 decomposition condition2}
involves five independent variables $r_{n+1}$, $r_{n}$,
$r_{n-1}$, $\gamma$ and the step-ratio limit $R_e$.
In general, we are not able to solve it exactly to
determine the optimal values of $\gamma$ so that the resulting step-ratio limit
$R_{e}$ is as large as possible.
As done in previous studies \cite{CalvoGrandeGrigorieff:1990,LiLiao:2021},
we consider a specific grid with constant step-ratio $r_k=r$ for a rough estimate of $\gamma$.
In such case the first condition \eqref{eq: BDF3 decomposition condition1} becomes
\begin{align*}
	\gamma\le\frac{-d_1(r,r)}{2\sqrt{r}d_2(r,r)}
	\quad \text{for all $r<R_{e}$.}
\end{align*}
An obvious choice  is
$\tilde{\gamma}=\frac{-d_1(\tilde{R}_{e},\tilde{R}_{e})}
{2\sqrt{\tilde{R}_{e}}d_2(\tilde{R}_{e},\tilde{R}_{e})}$ ($\tilde{\gamma}$
is decreasing with repect to $\tilde{R}_{e}$), but the parameter $\tilde{\gamma}$ is
not sufficient to ensure \eqref{eq: BDF3 decomposition condition1}.
Actually, we take $r_{n-1}=0$ and $r_{n+1}=r_{n}=\tilde{R}_e$ and get
\begin{align*}
	-\frac{\sqrt{\tilde{R}_e}}{\tilde{\gamma}}d_{1}(\tilde{R}_e,0)
	-\tilde{R}_ed_2(\tilde{R}_e,\tilde{R}_e)
	=\frac{-\tilde{R}_{e}^4(\tilde{R}_{e}^3+\tilde{R}_{e}^2-1)}
	{(\tilde{R}_{e}^2+\tilde{R}_{e}+1)(\tilde{R}_{e}^3
		+3 \tilde{R}_{e}^2+2 \tilde{R}_{e}+1)}
	<0
\end{align*}
if the step-ratio limit $\tilde{R}_e\ge1$, so that the desired condition \eqref{eq: BDF3 decomposition condition1}
is invalid. In turn, taking $r_{n-1}=0$ and $r_{n+1}=r_{n}=r$ in \eqref{eq: BDF3 decomposition condition1} yields
\begin{align*}
	\gamma\le\frac{-d_1(r,0)}{\sqrt{r}d_2(r,r)}
	\quad \text{for all $r<R_{e}$.}
\end{align*}
It introduces a reliable choice ($\bar{\gamma}$ is also decreasing with respect to $\bar{R}_{e}$)
\begin{align*}
		\bar{\gamma}:=\frac{-d_1(\bar{R}_{e},0)}{\sqrt{\bar{R}_{e}}d_2(\bar{R}_{e},\bar{R}_{e})}.
\end{align*}
Consider the second condition \eqref{eq: BDF3 decomposition condition2} with $\gamma=\bar{\gamma}$ and  $r_{n+1}=r_{n}=r_{n-1}=\bar{R}_e$. Thus we can determine the value of $\bar{R}_{e}$ by assuming that
\begin{align*}
	2d_0(\bar{R}_e,\bar{R}_e)-R_ed_2(\bar{R}_e,\bar{R}_e)
	+\bar{\gamma}^2\braB{\frac1{\bar{\gamma}}\sqrt{\bar{R}_e}d_1(\bar{R}_e,\bar{R}_e)
		+R_ed_2(\bar{R}_e,\bar{R}_e)}\\	
	+\frac1{\bar{\gamma}}\sqrt{\bar{R}_e}d_1(\bar{R}_e,\bar{R}_e)
	+\bar{R}_ed_2(\bar{R}_e,\bar{R}_e)=0.
\end{align*}
We solve this equation numerically and find the unique positive roots $\bar{R}_{e}\approx1.4965$
with the corresponding parameter $\bar{\gamma}\approx0.6924$.

To simplify the subsequent mathematical derivations, we fix the parameter $\gamma=7/10$
which is very close to $\bar{\gamma}\approx0.6924$.
The corresponding maximum step-ratio $R_{e}$ is determined by
$$d_1(R_{e},0)+\frac{7}{10}\sqrt{R_{e}}d_2(R_{e},R_{e})=0\quad
\text{or}\quad
\frac{10}{7\left(R_e+1\right)}-\frac{R_e^2\sqrt{R_e}}{R_e^2+R_e+1}=0.$$
We solve this equation numerically and find the unique positive root $R_e\approx1.4877.$
This choice is mainly because the restirction \eqref{eq: BDF3 decomposition condition1}
 should be necessary and sharp, while the inequality
  \eqref{eq: BDF3 decomposition condition2} can be relaxed appropriately.

We are now ready to prove the following lemma, which gives the discrete gradient structre.
Note that, the complex conditions
\eqref{eq: BDF3 decomposition condition1}-\eqref{eq: BDF3 decomposition condition2}
make the proof tedious and lengthy
and some technical lemmas are included in Appendix \ref{appendix: Technical results BDF3-DGS}.

\begin{lemma}\label{lem: BDF3-DGS}
	Define the following functions
	\begin{align}
			d_*(x,y):=&\,-\frac{10}{7}\sqrt{x}d_1(x,y)
			-\sqrt{xy}d_2(x,y),\quad\label{def: an}\\
		    p(x,y,z):=&\,2d_{0}(y,z)-\sqrt{yz}d_2(y,z)
		    -\frac{49}{100}d_*(y,z)-d_*(x,y),\label{def: An}\\
		     q(x,y,z):=&\,d_*(y,z)-\sqrt{xy}d_2(x,y),
		     \qquad \text{for $0<x,y,z<R_e$} \label{def: Bn}
	\end{align}
	If the step-ratios $0<r_k<R_e$, there exist two
	nonnegative functionals $G$ and $F$ such that
	\begin{align}\label{eq: BDF3-DGS}
		2v_n\tau_n\sum_{j=3}^nd_{n-j}^{(n)}v_j
		=&\,G[v_n,v_{n-1}]-G[v_{n-1},v_{n-2}]
		+F[v_n,v_{n-1},v_{n-2}]\quad\text{for $n\ge 3$,}
	\end{align}
	where the Lyapunov-type functional
	 \begin{align*}
	 	G[v_n,v_{n-1}]:=&\,d_*(r_{n+1},r_{n})\tau_nv_n^2
	 	+\sqrt{r_{n+1}r_{n}}d_2(r_{n+1},r_{n})
	 	\brab{\tfrac{7}{10}\sqrt{\tau_n}v_n-\sqrt{\tau_{n-1}}v_{n-1}}^2,
	 \end{align*}
 and the remainder term
 \begin{align*}
 	F[v_n,v_{n-1},v_{n-2}]:=&\,p(r_{n+1},r_{n},r_{n-1})\tau_nv_n^2+q(r_{n+1},r_{n},r_{n-1}) (\tfrac{7}{10}\sqrt{\tau_{n}}v_{n}-\sqrt{\tau_{n-1}}v_{n-1})^2\\
&\, +\sqrt{r_{n}r_{n-1}}d_2(r_{n},r_{n-1})
\brab{\sqrt{\tau_{n}}v_{n}-\tfrac{7}{10}\sqrt{\tau_{n-1}}v_{n-1}+\sqrt{\tau_{n-2}}v_{n-2}}^2
\ge \frac{\tau_n}{50}v_n^2.
 \end{align*}
	\end{lemma}

\begin{proof}
	According to Lemma \ref{lem: appendix-DGS-condition1} with
	$x:=r_{n+1}$, $y:=r_{n}$ and $z:=r_{n-1}$,
	the condition 	\eqref{eq: BDF3 decomposition condition1} holds
	for $\gamma=7/10$, that is,
		\begin{align*}
			q_{n+1}=&\,-\frac{10}{7} \sqrt{r_n}d_1^{(n)}
			-\sqrt{r_{n}r_{n-1}}d_2^{(n)}
		-\sqrt{r_{n+1}r_n}d_2^{(n+1)}\\
		=&\,q(r_{n+1},r_{n},r_{n-1})\ge0,\quad n\ge3.
	\end{align*}
Lemma \ref{lem: appendix-DGS-condition1} also implies that
$a(r_n,r_{n-1})\ge b(r_{n+1},r_n)>0$.
Applying Lemma \ref{lem: appendix-DGS-condition2} with
$x:=r_{n+1}$, $y:=r_{n}$ and $z:=r_{n-1}$,
one has
\begin{align*}
	p_{n+1}=&\,2d_{0}^{(n)}+\frac7{10}\sqrt{r_n}d_1^{(n)}
	-\frac{51}{100}\sqrt{r_{n}r_{n-1}}d_2^{(n)}
	+\frac{10}{7} \sqrt{r_{n+1}}d_1^{(n+1)}
	+\sqrt{r_{n+1}r_n}d_2^{(n+1)}\\
	=&\,p(r_{n+1},r_{n},r_{n-1})>\frac1{50},\quad n\ge3.
\end{align*}
Obviously, the condition \eqref{eq: BDF3 decomposition condition2}
holds for $\gamma=7/10$.
They imply that the discrete gradient structure
\eqref{eq: BDF3 decomposition assumption} holds,
that is,
\begin{align}\label{lemproof: BDF3-DGS}
	2w_n\sum_{j=3}^n
	\tilde{d}_{n-j}^{(n)}w_j=\widetilde{G}[w_n,w_{n-1}]
	-\widetilde{G}[w_{n-1},w_{n-2}]
	+\widetilde{F}[w_n,w_{n-1},w_{n-2}]
\end{align}
for $n\ge3$, where
\begin{align*}
	\widetilde{G}=&\,\braB  {-\tfrac{10}{7}
		\sqrt{r_{n+1}}d_{1}^{(n+1)}
		-\sqrt{r_{n+1}r_n}d_{2}^{(n+1)}}w_n^2
	+\sqrt{r_{n+1}r_n}d_{2}^{(n+1)}(\tfrac{7}{10}w_n-w_{n-1})^2\\
	=&\,d_*(r_{n+1},r_{n})w_n^2
	+\sqrt{r_{n+1}r_{n}}d_2(r_{n+1},r_{n})(\tfrac{7}{10}w_n-w_{n-1})^2,
\end{align*}
and the remainder term
 \begin{align*}
	\widetilde{F}=&\,
	p_{n+1}w_n^2+q_{n+1}(\tfrac{7}{10} w_n-w_{n-1})^2
	+\sqrt{r_{n}r_{n-1}}d_{2}^{(n)}\brab{w_n-\tfrac{7}{10} w_{n-1}+w_{n-2}}^2
	\ge \frac1{50}w_n^2.
\end{align*}
The claimed result follows from \eqref{lemproof: BDF3-DGS} immediately
by taking $w_k:=\sqrt{\tau_k}v_k$ and
\begin{align*}
	&G[v_n,v_{n-1}]:=\widetilde{G}\kbra{\sqrt{\tau_n}v_n,\sqrt{\tau_{n-1}}v_{n-1}},\\
	&F[v_n,v_{n-1},v_{n-2}]
	:=\widetilde{F}\kbra{\sqrt{\tau_n}v_n,\sqrt{\tau_{n-1}}v_{n-1},\sqrt{\tau_{n-2}}v_{n-2}}.
\end{align*}
It completes the proof.
\end{proof}

 On the uniform mesh with $d_0=11/6$, $d_1=-7/6$ and $d_2=1/3$,
 Lemma \ref{lem: BDF3-DGS} gives
\begin{align*}
	2v_n\braB{\frac{11}{6}v_n-\frac{7}{6}v_{n-1}+\frac{1}{3}v_{n-2}}
	=&\,\frac{4}{3}v_{n}^2+\frac1{3}(\tfrac{7}{10}v_n-v_{n-1})^2
	-\kbra{\frac{4}{3}v_{n-1}^2+\frac1{3}(\tfrac{7}{10}v_{n-1}-v_{n-2})^2}\\
	&\,+\frac{101}{75}v_n^2+(\tfrac{7}{10}v_n-v_{n-1})^2
	+\frac{1}{3}(v_n-\tfrac{7}{10}v_{n-1}+v_{n-2})^2.
\end{align*}
This decomposition arrives at a smaller bound $\frac{101}{75}\approx1.347$
than the optimal bound $\frac{96}{48}\approx1.979$ in \cite[Lemma 2.4]{LiaoTangZhou:2021bdf345}
for the minimum eigenvalue of the associated quadratic form. Actually, a sharp estimate of
the minimum eigenvalue is not the main purpose here. As seen, the main goal of
Lemma \ref{lem: BDF3-DGS} is to make the step-ratio limit $R_e$ as large as possible
on the basis of realizing a discrete gradient structure \eqref{eq: BDF3 energy quantity}.
According to the numerical tests at the beginning of this section, the
step-ratio limit $R_e\approx1.4877$ would be nearly optimal for having
a discrete gradient structure \eqref{lemproof: BDF3-DGS}, which implies
the positive definiteness  of discrete kernels $\big\{\tau_nd_{n-k}^{(n)}\big\}_{k=3}^n$.

\begin{lemma}\label{lem: BDF3-positive definite}
	If the step-ratios $0<r_k<R_e\approx1.4877$,
	the discrete kernels $\big\{\tau_nd_{n-k}^{(n)}\big\}_{k=3}^n$ are positive definite.
	in the sense that
	\begin{align*}
		2\sum_{k=3}^n\xi_k\sum_{j=3}^k\tau_kd_{k-j}^{(k)}\xi_j\ge \frac{1}{50}\sum_{k=3}^n\tau_k\xi_k^2		
		\quad\text{for $n\ge 3$.}
	\end{align*}	
\end{lemma}
\begin{proof}
	By taking $v_1=v_2=0$ and $v_j=\xi_j$ for $3\le j\le n$ in \eqref{eq: BDF3-DGS}
	and summing the resulting equalities form $n=3$ to $n=m$,
	one obtains the claimed inequality by replacing $m$ by $n$.
\end{proof}

\subsection{Energy dissipation law}

We now prove the energy ($H^1$ seminorm) stability of BDF3 scheme \eqref{eq: time-discrete IBVP}
for the dissipative case.
This property would be practically important when the
variable-step BDF3 scheme is applied to the gradient flow problems,
cf. the discussions \cite{ChenWangYanZhang:2019,LiaoJiZhang:2020pfc,
	LiaoTangZhou:2020bdf2}
 on variable-step BDF2 method.

\begin{theorem}\label{thm: enery stability}
	Assume that $\kappa\le0$ and $f\equiv0$. If $0<r_k<R_e\approx1.4877$ for $k\ge2$,
	the BDF3 scheme \eqref{eq: time-discrete IBVP} is unconditionally energy stable in the sense that
	\begin{align}\label{eq: discrete energy law}
		E^n\le E^{n-1},\quad\text{for $n\ge3$},
	\end{align}
where the (modified) discrete energy $E^n$ is defined by
\begin{align}\label{eq: discrete energy}
	E^n:=\varepsilon \mynormb{\nabla u^n}^2+\myinnerb{-\kappa u^{n},u^n}
	+\myinnerb{1,G\kbrab{\partial_{\tau} u^n,\partial_{\tau} u^{n-1}}}
	\quad \text{for $n\ge2$}.
\end{align}
\end{theorem}
\begin{proof}
Making the inner product of \eqref{eq: time-discrete IBVP} with $2\tau_n\partial_{\tau}u^n$, one obtains
\begin{align*}
	2\myinnerb{D_3u^{n},\tau_n\partial_{\tau}u^n}
	+2\varepsilon\myinnerb{\nabla u^{n},\nabla u^n-\nabla u^{n-1}}
	+2\myinnerb{-\kappa u^{n},u^n-u^{n-1}}=0,\quad\text{for $n\ge3$}.
\end{align*}
Taking $v_j=\partial_{\tau} u^j$ in Lemma \ref{lem: BDF3-DGS} gives
\begin{align*}
	2\myinnerb{D_3u^{n},\tau_n\partial_{\tau} u^n}\ge \myinnerb{1,G\kbrab{\partial_{\tau} u^n,\partial_{\tau} u^{n-1}}}-\myinnerb{1,G\kbrab{\partial_{\tau} u^{n-1},\partial_{\tau} u^{n-2}}}\quad\text{for $n\ge3$}.
\end{align*}
With the help of the inequality $2a(a-b)\ge a^2-b^2$, it is easy to obtain that
\begin{align*}
	E^n-E^{n-1}\le0,\quad\text{for $n\ge3$}.
\end{align*}
The proof is complete.
\end{proof}

\section{Stability and convergence analysis}
\setcounter{equation}{0}

In this section, we shall show the stability and convergence analysis of the variable-step BDF3 scheme.

\subsection{Properties of DOC kernels}

We first present the following lemma that shows the DOC-type kernels
$\big\{\tau_n\vartheta_{n-k}^{(n)}\big\}_{k=3}^n$
are positive definite.

\begin{lemma}\label{lem: BDF3-DOC-positive definite}
	If the step-ratios $0<r_k<R_e\approx1.4877$,
	the discrete kernels $\big\{\tau_n\vartheta_{n-k}^{(n)}\big\}_{k=3}^n$ are positive definite
	in the sense that, for any nonzero sequences $\{\xi_k\}$,
	\begin{align*}
		2\sum_{k=3}^n\xi_k\sum_{j=3}^k\tau_k\vartheta_{k-j}^{(k)}\xi_j
		>0,\quad\text{$n\ge 3$.}
	\end{align*}	
\end{lemma}
\begin{proof}
	For any nonzero sequences $\{\xi_k\}$,
	let $\eta_k:=\sum_{j=3}^k\vartheta_{k-j}^{(k)}\xi_j$ for $k\ge3$.
	Multiplying both sides of this equality by the BDF3 kernels
	$d_{n-k}^{(n)}$ and summing the index $k$
	from $k=3$ to $n$, we apply the orthogonality identity
	\eqref{eq: mutual BDFk-DOC identity} to get
	\begin{align*}
		\sum_{k=3}^nd_{n-k}^{(n)}\eta_k
		=&\,\sum_{k=3}^nd_{n-k}^{(n)}\sum_{j=3}^k\vartheta_{k-j}^{(k)}\xi_j
		=\sum_{j=3}^n\xi_j\sum_{k=j}^nd_{n-k}^{(n)}\vartheta_{k-j}^{(k)}=\xi_n
		\quad\text{for $n\ge 3$,}
	\end{align*}
where the summation order has been exchanged in the second equality.
Since the sequence $\{\eta_k\}$ is also nonzero,
Lemma \ref{lem: BDF3-positive definite} gives
\begin{align*}
	2\sum_{k=3}^n\xi_k\sum_{j=3}^k\tau_k\vartheta_{k-j}^{(k)}\xi_j
	=2\sum_{k=3}^n\eta_k\tau_k\sum_{j=3}^kd_{k-j}^{(k)}\eta_j
	>0\quad\text{for $n\ge 3$.}
\end{align*}
The claimed inequality is verified.
\end{proof}

As noted in \cite{LiLiao:2021}, the DOC kernels
$\vartheta_{j}^{(n)}$ in \eqref{eq: BDFk-DOC procedure}
 are not always positive, but they decay exponentially
 such that the absolute summations of DOC kernels
 are uniformly bounded, as stated in the following lemma.

\begin{lemma}\cite[Lemma 3.1]{LiLiao:2021}\label{lem: BDF3 orthogonal formula}
	If the step ratios $r_k\le R_e$ for $k\ge2$, there exists a positive constant $K_3$
	such that the DOC kernels $\vartheta_{n-j}^{(n)}$ in \eqref{eq: BDFk-DOC procedure}
	satisfy
	\begin{align*}
		\sum_{j=3}^{n}\absb{\vartheta_{n-j}^{(n)}}\le K_3
		\quad\text{and}\quad
		\sum_{j=i}^{n}\absb{\vartheta_{j-i}^{(j)}}\le K_3\quad\text{for $n\ge 3\;(i\ge3)$,}
	\end{align*}
where $K_3$ is independent of the time $t_n$ and the step ratios $r_k\in(0,R_e]$ $(k\ge2)$.
\end{lemma}

\subsection{$L^2$ norm stability}
We are now ready to present the $L^2$ norm stability.
\begin{theorem}\label{thm: BDF3 L2 stability-dissipatve}
	If the step ratios $0<r_k<R_e$ for $k\ge2$, the BDF3 solution
	of \eqref{eq: time-discrete IBVP} with $\kappa<0$
	is mesh-robustly stable in the $L^2$ norm, that is,
	\begin{align*}
		\mynormb{u^n}\le&\, \mynormb{u^2}+K_3\tau\mynormb{\partial_{\tau}u^1}+4K_3\tau\mynormb{\partial_{\tau}u^2}
			+2\sum_{k=3}^{n}\tau_k\sum_{j=3}^{k}\absb{\theta_{k-j}^{(k)}}\mynormb{f^j}\\
			\le&\, \mynormb{u^2}+K_3\tau\mynormb{\partial_{\tau}u^1}+4K_3\tau\mynormb{\partial_{\tau}u^2}
			+2K_3t_n\max_{3\le k\le n}\mynormb{f^k}\quad\text{for $n\ge3$.}
	\end{align*}
\end{theorem}

\begin{proof}Thanks to Lemma \ref{lem: BDF3 orthogonal formula}, it remains to verify the first estimate.
	Making the inner product of the equation \eqref{eq: DOC-BDF3 equation}  with $2\tau_nu^n$,
	and summing the resulting equality from $n=3$ to $m$, one has
		\begin{align}\label{thmproof: BDF3 L2 stability-dissipatve-0}
		2\sum_{k=3}^{m}\tau_k\myinnerb{u^{k},\partial_{\tau}u^{k}}
		=&\,-2\sum_{k=3}^{m}\tau_k\myinnerb{u^{k},\mathcal{I}_3^k[u]}
		-2\varepsilon\sum_{k=3}^{m}
		\sum_{j=3}^{k}\tau_k\theta_{k-j}^{(k)}\myinnerb{\nabla u^{k},\nabla u^j}\nonumber\\
		&\,-2\sum_{k=3}^{m}\sum_{j=3}^{k}\tau_k\theta_{k-j}^{(k)}\myinnerb{u^{k},(-\kappa)u^j}
		+2\sum_{k=3}^{m}\sum_{j=3}^{k}\tau_k\myinnerb{u^{k},\theta_{k-j}^{(k)}f^j}
	\end{align}
for $m\ge3$. Lemma \ref{lem: BDF3-DOC-positive definite} leads to
	\begin{align*}
		2\sum_{k=3}^{m}\tau_k\myinnerb{u^{k},\partial_{\tau}u^{k}}
		\le-2\sum_{k=3}^{m}\tau_k\myinnerb{u^{k},\mathcal{I}_3^k[u]}
		+2\sum_{k=3}^{m}\sum_{j=3}^{k}\tau_k\myinnerb{u^{k},\theta_{k-j}^{(k)}f^j}
		\quad\text{for $m\ge3$.}
	\end{align*}
	Note that, $2\tau_k\myinnerb{u^{k},\partial_{\tau}u^{k}}\ge\mynormt{u^{k}}^2-\mynormt{u^{k-1}}^2$.
	Then applying the Cauchy-Schwarz inequality and Lemma \ref{lem: BDF3 orthogonal formula}, we get	
	\begin{align*}
		\mynormb{u^{m}}^2\le  \mynormb{u^{2}}^2
		+2\sum_{k=3}^{m}\tau_k\mynormb{u^{k}}\mynormb{\mathcal{I}_3^k[u]}
		+2\sum_{k=3}^{m}\sum_{j=3}^{k}\tau_k\absb{\theta_{k-j}^{(k)}}\mynormb{u^{k}}\mynormb{f^j}
		\quad\text{for $m\ge3$.}
	\end{align*}
	Taking some integer $m_0$ ($2\le m_0\le m$) such that $\mynormb{u^{m_0}}=\max_{2\le k\le m}\mynormb{u^{k}}$. Taking $m:=m_0$ in the above inequality, one gets
	\begin{align*}
		\mynormb{u^{m_0}}^2\le  \mynormb{u^{2}}\mynormb{u^{m_0}}
		+2\mynormb{u^{m_0}}\sum_{k=3}^{m_0}\tau_k\mynormb{\mathcal{I}_3^k[u]}
		+2\mynormb{u^{m_0}}\sum_{k=3}^{m_0}\sum_{j=3}^{k}\tau_k\absb{\theta_{k-j}^{(k)}}\mynormb{f^j}
	\end{align*}
	and thus
	\begin{align}\label{thmproof: BDF3 L2 stability-dissipatve-1}
		\mynormb{u^{m}}\le&\,\mynormb{u^{m_0}}\le  \mynormb{u^{2}}
		+2\tau\sum_{k=3}^{m}\mynormb{\mathcal{I}_3^k[u]}
		+2\sum_{k=3}^{m}\sum_{j=3}^{k}\tau_k\absb{\theta_{k-j}^{(k)}}\mynormb{f^j}
		\quad\text{for $m\ge3$.}
	\end{align}

Now we evaluate the term $\sum_{k=3}^m\absb{\mathcal{I}_{3}^k[u]}$
stemmed from the starting values.
Recalling the definition \eqref{eq: BDF3 kernels} of discrete BDF3 kernels
 with the increasing property \eqref{eq: d0d1d2-increasing} of $\abs{d_1(x,y)}$ and $d_2(x,y)$,
it is easy to check that
$$d^{(3)}_{2}\le d_2(R_e,R_e)\le \frac1{2}\quad\text{and}\quad
\absb{d^{(3)}_{1}}+d^{(4)}_{2}\le -d_1(R_e,R_e)+d_2(R_e,R_e)\le 2.$$
Thus we apply the formula \eqref{eq: initial effect BDF3-DOC}
and  Lemma \ref{lem: BDF3 orthogonal formula} to get
\begin{align}\label{thmproof: BDF3 L2 stability-initial}
	\sum_{k=3}^m\absb{\mathcal{I}_{3}^k[u]}\le&\,
	\absb{\partial_{\tau}u^2}
	\sum_{k=3}^m\sum_{i=3}^k\absb{\vartheta_{k-i}^{(k)}}\absb{d^{(i)}_{i-2}}
	+\absb{\partial_{\tau}u^1}
	\sum_{k=3}^m\absb{\vartheta_{k-3}^{(k)}}d^{(3)}_{2}\nonumber\\
	=&\,\absb{\partial_{\tau}u^2}
	\sum_{i=3}^m\absb{d^{(i)}_{i-2}}\sum_{k=i}^m\absb{\vartheta_{k-i}^{(k)}}
	+\absb{\partial_{\tau}u^1}
	d^{(3)}_{2}\sum_{k=3}^m\absb{\vartheta_{k-3}^{(k)}}\nonumber\\
	\le&\, 2K_3\absb{\partial_{\tau}u^2}+\frac1{2}K_3\absb{\partial_{\tau}u^1}\quad\text{for $m\ge3$.}
\end{align}
Inserting this estimate \eqref{thmproof: BDF3 L2 stability-initial}
into \eqref{thmproof: BDF3 L2 stability-dissipatve-1},
one obtains the first result and completes the proof.
\end{proof}

\begin{theorem}\label{thm: BDF3 L2 stability}
Let $\kappa$ be bounded such that
$\abs{\kappa}\le \kappa^*$. If the step ratios $0<r_k<R_e$ $(k\ge2)$ with
the maximum time-step $\tau\le 1/(4K_3\kappa^*)$, the BDF3 solution
 of \eqref{eq: time-discrete IBVP} satisfies
 \begin{align*}
\mynormb{u^n}\le 2\exp(4K_3\kappa^*t_{n-1})\braB{ \mynormb{u^2}+K_3\tau\mynormb{\partial_{\tau}u^1}+4K_3\tau\mynormb{\partial_{\tau}u^2}
	+2K_3t_n\max_{3\le k\le n}\mynormb{f^k}}
\end{align*}
for $3\le n\le N$. Thus the BDF3 scheme is mesh-robustly stable in the $L^2$ norm.
\end{theorem}

\begin{proof} This proof is an easy extension to that of
	Theorem \ref{thm: BDF3 L2 stability-dissipatve}.
	By following from the equality \eqref{thmproof: BDF3 L2 stability-dissipatve-0},
	it is not difficult to get
		\begin{align*}
		\mynormb{u^{m}}^2\le&\,  \mynormb{u^{2}}^2
		+2\sum_{k=3}^{m}\tau_k\mynormb{u^{k}}\mynormb{\mathcal{I}_3^k[u]}
		+2\kappa^{*}\sum_{k=3}^{m}\sum_{j=3}^{k}\tau_k\absb{\theta_{k-j}^{(k)}}\mynormb{u^{k}}\mynormb{u^{j}}\\
		&\,+2\sum_{k=3}^{m}\sum_{j=3}^{k}\tau_k\absb{\theta_{k-j}^{(k)}}\mynormb{u^{k}}\mynormb{f^j}
		\quad\text{for $m\ge3$.}
	\end{align*}
Taking some integer $m_1$ ($2\le m_1\le m$) such that
$\mynormb{u^{m_1}}=\max_{2\le k\le m}\mynormb{u^{k}}$.
 Taking $m:=m_1$ in the above inequality, one gets
\begin{align*}
	\mynormb{u^{m_1}}\le  \mynormb{u^{2}}
	+2\sum_{k=3}^{m_1}\tau_k\mynormb{\mathcal{I}_3^k[u]}
	+2\kappa^{*}\sum_{k=3}^{m_1}\tau_k\mynormb{u^{k}}\sum_{j=3}^{k}\absb{\theta_{k-j}^{(k)}}
	+2\sum_{k=3}^{m_1}\sum_{j=3}^{k}\tau_k\absb{\theta_{k-j}^{(k)}}\mynormb{f^j}.
\end{align*}
Applying Lemma \ref{lem: BDF3 orthogonal formula} and the initial estimate
\eqref{thmproof: BDF3 L2 stability-initial}, we can derive that
\begin{align*}
	\mynormb{u^{m}}\le&\,\mynormb{u^{2}}
	+K_3\tau\mynormb{\partial_{\tau}u^1}+4K_3\tau\mynormb{\partial_{\tau}u^2}	
	+2K_3\kappa^{*}\sum_{k=3}^{m}\tau_k\mynormb{u^{k}}
	+2K_3t_n\max_{3\le k\le m}\mynormb{f^k}
\end{align*}
for $3\le m\le N$. With the time-step condition $\tau\le 1/(4K_3\kappa^{*})$, it arrives at
\begin{align*}
	\mynormb{u^{n}}\le&\,2\mynormb{u^{2}}
	+2K_3\tau\mynormb{\partial_{\tau}u^1}+8K_3\tau\mynormb{\partial_{\tau}u^2}	
	+4K_3\kappa^{*}\sum_{k=3}^{n-1}\tau_k\mynormb{u^{k}}
	+4K_3t_n\max_{3\le k\le n}\mynormb{f^k}
\end{align*}
Then the standard Gr\"{o}nwall inequality completes the proof.
\end{proof}

The above theorems remove the unbounded quantity
$\Gamma_n$ in \cite{CalvoGrigorieff:2002} completely.
They show that the variable-step BDF3 scheme is surprisingly  stable
with respect to the changes of time-steps if the step ratios
satisfy a sufficient condition $r_k<R_e$.
Extensive tests on random time meshes in Section 4 suggest that
this step-ratio constraint is far from necessary for the stability.

\subsection{$L^2$ norm convergence}

We finally present the $L^2$ norm convergence. For the auxiliary functions \eqref{funs: d0}-\eqref{funs: d2}, it is easy  to check that
\begin{align}
	&d_{0}(x,y)+d_{1}(x,y)+d_{2}(x,y)=1,\label{funs: d012-sum1}\\
	&d_{0}(x,y)+\frac{1+2x}{x}d_{1}(x,y)+\frac{1+2y+2xy}{xy}d_{2}(x,y)=0,\label{funs: d012-sum2}\\
	&d_{0}(x,y)+\frac{1+3x+3x^2}{x^2}d_{1}(x,y)
	+\frac{1+3y+3xy+3y^2+6xy^2+3x^2y^2}{x^2y^2}d_{2}(x,y)=0.\label{funs: d012-sum3}
\end{align}
Consider a smooth function $v$ and let $\zeta^j[v]:=D_3v(t_j)-v'(t_j)$ be the truncation error at $t_j$ $(j\ge3)$
of variable-step BDF3 formula. By applying the Taylor's expansion with integral remainder,
one can apply the identies \eqref{funs: d012-sum1}-\eqref{funs: d012-sum3} to find that
\begin{align}
	\zeta^j[v]
	=&\,\sum_{i=j-2}^j\frac{1}{6\tau_{i}}
	\int_{t_{i-1}}^{t_{i}}\mathcal{K}_{j,j-i}(t)v^{(4)}(t)\zd{t},
\end{align}
where the involoved integral kernels read
\begin{align*}
	\mathcal{K}_{j,0}(t)=&\,\brab{d_0^{(j)}-r_jd_1^{(j)}}(t_{j-1}-t)^3
	+r_j\brab{d_1^{(j)}-r_{j-1}d_2^{(j)}}(t_{j-2}-t)^3
	+r_jr_{j-1}d_2^{(j)}(t_{j-3}-t)^3,\\
	\mathcal{K}_{j,1}(t)=&\,\brab{d_1^{(j)}-r_{j-1}d_2^{(j)}}(t_{j-2}-t)^3+r_{j-1}d_2^{(j)}(t_{j-3}-t)^3,\\
		\mathcal{K}_{j,2}(t)=&\,d_2^{(j)}(t_{j-3}-t)^3.
\end{align*}
Reminding the increasing property \eqref{eq: d0d1d2-increasing}, it is not difficult to prove that
there exists a bounded constant $K_v$ such that
\begin{align}\label{eq: BDF3 truncation error bound}	
	\absb{\zeta^j[v]}\le K_v\tau^3	\quad	\text{for $3\le j\le N$},
\end{align}
where $K_v$ is always dependent on the function $v$, but independent of the time $t_n$,
the step sizes $\tau_n$ and the step ratios $r_n$ (even when $r_n$ approaches the limit $R_e$).

Let $\tilde{u}^n:=u(t_n)-u^n$ be the solution error of
the variable-step BDF3 scheme \eqref{eq: time-discrete IBVP}.
We have the error equation
\begin{align}\label{eq: BDF3 error equation}	
	D_3\tilde{u}^{k}=\varepsilon\Delta \tilde{u}^{k}
	+\kappa \tilde{u}^{k}+\zeta^k[u],\quad\text{for $1\le k\le N$}
\end{align}
together with the initial conditions~$\tilde{u}^0$, $\tilde{u}^1$ and $\tilde{u}^2$.
Then applying the priori stability estimate in Theorem \ref{thm: BDF3 L2 stability}
to the error equation \eqref{eq: BDF3 error equation}, we can use
the error bound \eqref{eq: BDF3 truncation error bound} to
verify the following convergence result.

\begin{theorem}\label{thm: BDF3 L2 convergence}
	Assume that the solution of \eqref{eq: diffusion problem} is smooth enough in time.
	If the time-step ratios $0<r_k<R_e$ $(k\ge2)$ with
	the maximum time-step $\tau\le 1/(4K_3\kappa^*)$,
	the solution error $\tilde{u}^n=u(t_n)-u^n$ of
	the variable-step BDF3 scheme \eqref{eq: time-discrete IBVP} satisfies
	\begin{align*}
		\mynormb{\tilde{u}^n}\le 2\exp(4K_3\kappa^*t_{n-1})\braB{ \mynormb{\tilde{u}^2}+K_3\tau\mynormb{\partial_{\tau}\tilde{u}^1}
			+4K_3\tau\mynormb{\partial_{\tau}\tilde{u}^2}
			+2K_3K_ut_n\tau^3}
	\end{align*}
	for $3\le n\le N$. Here the constants $K_3$ and $K_u$ are independent of the time $t_n$,
	the step sizes $\tau_n$ and the step ratios $r_n$ (even when $r_n$ approaches the limit $R_e$).
	Thus the BDF3 scheme is mesh-robustly convergent in the $L^2$ norm.
\end{theorem}

To start the third-order stiff solver, one can apply a
third-order Runge-Kutta method to compute the starting solutions $u^1$ and $u^2$.
Our error estimate in Theorem \ref{thm: BDF3 L2 convergence} also
implies that a second-order starting scheme for computing $u^1$ and $u^2$
would be adequate to achieve the overall third-order accuracy
since it can generate third-order accurate solutions
at the first two levels.

\section{Numerical experiments}
\setcounter{equation}{0}

We shall present in this section some numerical examples to support our theoretical findings. To this end,
we consider the heat equation $\partial_tu-\varepsilon\Delta u=f$
on the square domain $\Omega=(0,2\pi)^2$ with periodic boundary conditions.
We choose the exterior force $f$ and the diffusive coefficient $\varepsilon=0.1$
such that the equation yields a smooth solution $u=\cos(t)\sin(x)\sin(y)$.

To start the three-step stiff solver, we use the variable-step BDF2 method
and a two-stage third-order singly diagonally implicit Runge-Kutta method in our numerical implementations.
The numerical stability and convergence are tested until time $T=1$ in two scenarios:
\begin{enumerate}
	\item[(a)] The periodic time steps  $\{\tau_1,\mu\tau_1,\cdots,\tau_1,\mu\tau_1,\cdots,\tau_1,\mu\tau_1\}$
	with a constant $\mu>1$, where $\tau_1=2/(N(1+\mu))$
	and the maximum step-ratio $r_{\max}=r_{2j}=\mu$ $(j=1,2,\cdots,N/2)$.	
	\item[(b)] The random time steps $\tau_k=\epsilon_k/\sum_{k=1}^N\epsilon_k$,
	where $\epsilon_k\in (0,1)$ are uniformly distributed random numbers.
\end{enumerate}

\begin{table}[htb!]	
	\begin{center}
		\caption{The BDF3 solutions on periodic meshes with $\mu=2R_e$
			starting by a third-order RK}
			\label{table:BDF3-RK periodic 2Re}\vspace*{0.3pt}
		\def\temptablewidth{0.6\textwidth}
		{\rule{\temptablewidth}{0.5pt}}
		\begin{tabular*}{\temptablewidth}{@{\extracolsep{\fill}}ccccc}
			$N$   &$\tau(N)$   &$e(N)$     &Order  &$N_1$\\
			\midrule
			80    &1.87e-02	 &1.12e-06	&-- 	 	 &40\\
			160   &9.36e-03	 &1.42e-07	&2.98	 	 &80\\
			320	  &4.68e-03	 &1.78e-08	&2.99	 	 &160\\	
			640	  &2.34e-03	 &2.23e-09	&3.00		 &320\\
			1280  &1.17e-03	 &2.80e-10	&3.00		 &640
		\end{tabular*}
		{\rule{\temptablewidth}{0.5pt}}
	\end{center}
\end{table}	

\begin{table}[htb!]	
	\begin{center}
		\caption{The BDF3 solutions on periodic meshes with $\mu=4R_e$
			starting by the BDF2 method}
		\label{table:BDF3-BDF2 periodic 4Re}
		\vspace*{0.3pt}
		\def\temptablewidth{0.6\textwidth}
		{\rule{\temptablewidth}{0.5pt}}
		\begin{tabular*}{\temptablewidth}{@{\extracolsep{\fill}}ccccc}
			$N$   &$\tau(N)$   &$e(N)$     &Order  &$N_1$\\
			\midrule
			80    &2.14e-02	   &1.10e-06   &--  	&40\\
			160   &1.07e-02	   &1.39e-07   &2.98	&80\\
			320	  &5.35e-03	   &1.74e-08   &2.99    &160\\	
			640	  &2.68e-03	   &2.19e-09   &3.00    &320\\
			1280  &1.34e-03	   &2.74e-10   &2.99    &640
		\end{tabular*}
		{\rule{\temptablewidth}{0.5pt}}
	\end{center}
\end{table}	

\begin{table}[htb!]	
	\begin{center}
			\caption{The BDF3 solutions on random meshes starting by a third-order RK}
		\label{table:BDF3-RK random}
		\vspace*{0.3pt}
		\def\temptablewidth{0.7\textwidth}
		{\rule{\temptablewidth}{0.5pt}}
		\begin{tabular*}{\temptablewidth}{@{\extracolsep{\fill}}cccccc}
			$N$   &$\tau(N)$    &$e(N)$     &Order  &$r_{\max}$ &$N_1$\\
			\midrule
			80    &2.51e-02	   &1.61e-06    &-- 	&28.32	    &29\\
			160   &1.29e-02	   &2.18e-07	&2.88	&167.21     &53\\
			320	  &6.32e-03	   &2.75e-08	&2.99	&401.76     &110\\	
			640	  &3.18e-03	   &3.53e-09	&2.96	&1656.74    &206\\
			1280  &1.53e-03	   &4.36e-10	&3.02	&1584.01    &420
		\end{tabular*}
		{\rule{\temptablewidth}{0.5pt}}
	\end{center}
\end{table}

\begin{table}[htb!]
	\begin{center}
		\caption{The BDF3 solutions on random meshes starting by the BDF2 method}
		\label{table:BDF3-BDF2 random}\vspace*{0.3pt}
		\def\temptablewidth{0.7\textwidth}
		{\rule{\temptablewidth}{0.5pt}}
		\begin{tabular*}{\temptablewidth}{@{\extracolsep{\fill}}cccccc}
			$N$   &$\tau(N)$   &$e(N)$     &Order  &$r_{\max}$ &$N_1$\\
			\midrule
			80    &2.42e-02	   &1.71e-06	&--	&746.55	   &13\\
			160   &1.22e-02	   &2.42e-07    &2.82	&110.90	   &26\\
			320	  &6.54e-03	   &2.90e-08	&3.06	&79.85	   &70\\	
			640	  &3.18e-03	   &3.31e-09	&3.13	&371.22	   &125\\
			1280  &1.57e-03	   &4.44e-10	&2.90	&1321.80   &254
		\end{tabular*}
		{\rule{\temptablewidth}{0.5pt}}
	\end{center}
\end{table}	

We record the $L^2$ norm error
$e(N):=\max_{1\leq{n}\leq{N}}\|v(t_n)-v^n\|$
in each run and compute the numerical order of convergence by
$$\text{Order}\approx\frac{\log\bra{e(N)/e(2N)}}{\log\bra{\tau(N)/\tau(2N)}}$$
where $\tau(N)$ denotes the maximum time-step size for total $N$ subintervals.

Numerical results on the periodic time meshes
are listed in Tables
\ref{table:BDF3-RK periodic 2Re}-\ref{table:BDF3-BDF2 periodic 4Re},
in which we also record the number (denote by $N_1$)
of time levels with the step ratios $r_k\ge R_e\approx 1.4877$.
We observe that (i) the numerical solution is stable
even if there are 50\% of step-ratios greater than our theoretical restriction;
(ii) both a third-order SDIRK method and
the second-order BDF2 method are enough to achieve the third-order accuracy,
as predicted in Theorem \ref{thm: BDF3 L2 convergence}.
Table \ref{table:BDF3-RK random}-\ref{table:BDF3-BDF2 random}
record the numerical results on random time meshes.
We see that variable-step BDF3 method is mesh-robust with a  desired convergence rate,
even if many of step-ratios are much greater than our theoretical limit, and this well be further investigated in our future studies.


\section*{Acknowledgements}
The authors would like to thank Dr. Ji Bingquan
and Dr. Wang Jindi for their help on extensive numerical tests
on random time meshes.

\appendix
\section{Technical results for Lemma \ref{lem: BDF3-DGS}}
\label{appendix: Technical results BDF3-DGS}
\setcounter{equation}{0}

By using the definitions \eqref{funs: d0}-\eqref{funs: d2}, it is easy to check that
\begin{align}\label{eq: d0d1d2-increasing}
	&\frac{\partial }{\partial x}\abs{d_{\nu}(x,y)}
	>0\quad\text{and}\quad \frac{\partial}{\partial y}\abs{d_{\nu}(x,y)}
	>0\quad\text{for $\nu=0,1,2$ and $x,y>0.$}
\end{align}
That is, the functions $d_0(x,y)$, $-d_1(x,y)$ and $d_2(x,y)$ are increasing with respect to $x,y>0$.


This appendix presents some technical lemmas related to the functions
$d_{0}$, $d_{1}$ and $d_{2}$ defined by \eqref{funs: d0}-\eqref{funs: d2}, respectively.
The subsequent analysis is somewhat technically complex and
the mathematical derivations have been checked carefully by a symbolic calculation software.

\begin{lemma}\label{lem: appendix-DGS-condition1}
	For the function $q$ defined by \eqref{def: An}, it holds that
	\begin{align*}
		q=-\frac{10}{7} \sqrt{y}d_1(y,z)-\sqrt{yz}d_2(y,z)-\sqrt{xy}d_2(x,y)>0\quad
		\text{for $0<x,y,z<R_e$.}
	\end{align*}
\end{lemma}
\begin{proof}
	We consider the auxiliary function
	\begin{align*}
		\eta(x,y,z):=&\,\frac1{y^3}\kbra{-\frac{10}{7} yd_1(y^2,z^2)-yzd_2(y^2,z^2)-xyd_2(x^2,y^2)}\\
		=&\,\frac{10}{7\left(y^2+1\right)}
		-\frac{x^3 \left(x^2+1\right) y^2}{\left(y^2+1\right) \left(x^2y^2+y^2+1\right)}\\
		&\,+\frac{10 z^2 \left(y^2z^2+2z^2+1\right)}{7\left(z^2+1\right) \left(y^2z^2+z^2 +1\right)}-\frac{\left(y^2+1\right) z^5}{\left(z^2+1\right) \left(y^2z^2+z^2 +1\right)}.
	\end{align*}
It remains to verify $\eta(x,y,z)>0$ for $0< x, y,z < \sqrt{R_{e}}$. Simple calculations give
	\begin{align}\label{lemproof: appendix-eta-x}
		\frac{\partial \eta}{\partial x}=-\frac{x^2 y^2 \left[3 x^4 y^2+x^2 \left(6 y^2+5\right)+3 \left(y^2+1\right)\right]}{\left(y^2+1\right) \left(x^2y^2+y^2+1\right)^2}<0
	\end{align}
	such that $\eta(x,y,z)>\eta(\sqrt{R_{e}},y,z)$ for $0< x, y,z < \sqrt{R_{e}}$.
	Moreover, we get
		\begin{align}\label{lemproof: appendix-DGS-condition1-1}
		\frac{\partial \eta}{\partial y}=-\frac{2y}{7(y^2+1)^2}\kbrab{-\eta_1(x,y)+\eta_2(y,z)},
	\end{align}
	where $\eta_1$ and $\eta_2$ are defined by
	\begin{align*}
		\eta_1(x,y):=&\,\frac{7 x^3 \left(x^2+1\right) \left(\left(x^2+1\right) y^4-1\right)}{\left(x^2y^2+y^2+1\right)^2},\\
		\eta_2(y,z):=&\,\frac{20 \left(y^2+1\right)^2 z^6+7 \left(y^2+1\right)^2 z^5+10 \left(y^4+4 y^2+3\right) z^4+10 \left(2 y^2+3\right) z^2+10}{\left(z^2+1\right) \left(y^2z^2+z^2+1\right)^2}.
	\end{align*}
For the function $\eta_1$, one has
\begin{align*}
	\frac{\partial \eta_1}{\partial y}=
	\frac{28 x^3 \left(x^2+1\right)^2 \left(y^3+y\right)}{\left(x^2y^2+y^2+1\right)^3}>0
	\quad\text{for $0<x,y< \sqrt{R_{e}}$}
\end{align*}
and
\begin{align*}
	\frac{\partial \eta_1(x,\sqrt{R_{e}}) }{\partial x}
=&\,\frac{7x^2}{\left(R_ex^2+R_e+1\right)^3}\left[\Big.\Big.3R_e^3 x^6 +\left(9 R_e^2+7 R_e-1\right) R_ex^4 \right.\\
&\,\left.+
\left(9 R_e^3+10 R_e^2-4R_e-5\right)x^2 +3\left(R_e-1\right) \left(R_e+1\right)^2\right]>0
\end{align*}
such that
\begin{align}\label{lemproof: appendix-DGS-condition1-2}
	\eta_1(x,y)<\eta_1(x,\sqrt{R_{e}})<\eta_1(\sqrt{R_{e}},\sqrt{R_{e}})<7\quad
	\text{for $0< x, y< \sqrt{R_{e}}$}.
\end{align}
For the function $\eta_2$, one has
\begin{align*}
	\frac{\partial \eta_2}{\partial y}=
	\frac{4 y \left(y^2+1\right) z^5 (10 z+7)}{\left(z^2+1\right) \left(y^2z^2+z^2+1\right)^3}>0
	\quad\text{for $0<y,z< \sqrt{R_{e}}$}
\end{align*}
and
\begin{align*}
	\frac{\partial \eta_2(0,z)}{\partial z} =
	\frac{z^4(7z-30+\sqrt{1145})(30+\sqrt{1145}-7z)}{7\left(z^2+1\right)^4}>0
	\quad\text{for $0<z< \sqrt{R_{e}}$}
\end{align*}
such that
$$\eta_2(y,z)>\eta_2(0,z)>\eta_2(0,0)=10\quad
\text{for $0< y, z< \sqrt{R_{e}}$}.$$
Thanks to \eqref{lemproof: appendix-DGS-condition1-2},  one has
$$\eta_1(x,y)<\eta_2(y,z)\quad\text{for $0<x, y, z< \sqrt{R_{e}}$}.$$
Thus the fromula \eqref{lemproof: appendix-DGS-condition1-1} shows that
$\frac{\partial \eta}{\partial y}<0$ for $0< x, y,z < \sqrt{R_{e}}$,
see Figure 1 (a), and then
\begin{align}\label{lemproof: appendix-DGS-condition1-3}
\eta(x,y,z)>\eta(\sqrt{R_{e}},y,z)>\eta(\sqrt{R_{e}},\sqrt{R_{e}},z)\quad
\text{for $0< x,y, z< \sqrt{R_{e}}$}.
\end{align}
It remains to examine the following function with respect to $0<z< \sqrt{R_{e}}$,
\begin{align*}
	\eta(\sqrt{R_{e}},\sqrt{R_{e}},z)
	=\frac{10}{7\left(R_e+1\right)}-\frac{R_e^2\sqrt{R_e}}{R_e^2+R_e+1}
	-\frac{7 \left(R_e+1\right) z^5-10 \left(R_e+2\right) z^4-10 z^2}
	{7 \left(z^2+1\right) \left(\left(R_e+1\right) z^2+1\right)}.
\end{align*}
We solve the equation $\frac{\zd}{\zd z}\eta(\sqrt{R_{e}},\sqrt{R_{e}},z)=0$
numerically and find two real solutions $z_0^*=0$ and $z_1^*=1.07229$. Thus
\begin{align*}
	\eta(\sqrt{R_{e}},\sqrt{R_{e}},z)
	>&\,\min\left\{\eta(\sqrt{R_{e}},\sqrt{R_{e}},0),\eta(\sqrt{R_{e}},\sqrt{R_{e}},z_1^*),
	\eta(\sqrt{R_{e}},\sqrt{R_{e}},\sqrt{R_{e}})\right\}\\
	=&\,\eta(\sqrt{R_{e}},\sqrt{R_{e}},0)\approx 0.0000017>0.
\end{align*}
Thus the fact \eqref{lemproof: appendix-DGS-condition1-3} arrives at the claimed inequality
and completes the proof.
\end{proof}

\begin{figure}[htb!]	
	\begin{center}
\includegraphics[width=2.5in]{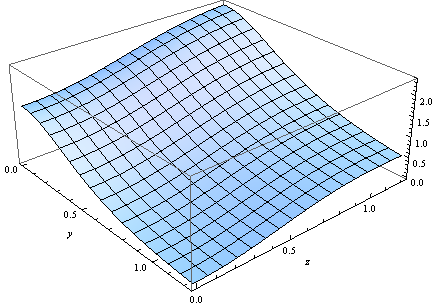}\quad
\includegraphics[width=2.5in]{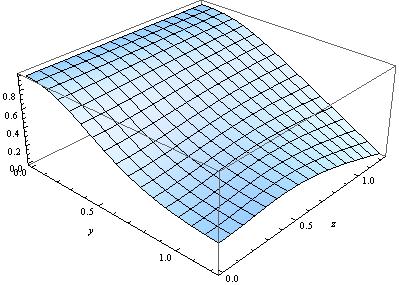}
	\end{center}
\caption{Surfaces of $\eta(\sqrt{R_e},y,z)$ and $\zeta(\sqrt{R_e},y,z)$ over the domain $(0,\sqrt{R_{e}})^2$.}
\label{fig: two surfaces}
\end{figure}

\begin{lemma}\label{lem: appendix-DGS-condition2}
		For the function $p$ defined by \eqref{def: Bn} for $0<x,y,z<R_{e}$, it holds that
	\begin{align*}
		p=&\,2d_0(y,z)+\frac{7}{10}\sqrt{y}d_1(y,z)-\frac{51}{100}\sqrt{yz}d_2(y,z)
		+\frac{10}{7}\sqrt{x}d_1(x,y)+\sqrt{xy}d_2(x,y)>\frac1{50}.
	\end{align*}
\end{lemma}
\begin{proof}Reminding the function $d_*$ defined by \eqref{def: an},
	we consider an auxiliary function
	$$\zeta(x,y,z):=\frac{\zeta_1(y,z)+\zeta_2(x,y)}{1+y^2}\quad\text{for $0<x,y,z< \sqrt{R_{e}}$},$$
	where
\begin{align*}	
	\zeta_1(y,z):=&\,2d_0(y^2,z^2)+\frac{7}{10}y d_1(y^2,z^2)-\frac{51}{100}yzd_2(y^2,z^2),\\
	\zeta_2(x,y):=&\,-d_*(x^2,y^2)=\frac{10}{7}xd_1(x^2,y^2)+xyd_2(x^2,y^2).
\end{align*}
Obviously, one has
	\begin{align}\label{lemproof: appendix-DGS-condition2-0}
	p(x^2,y^2,z^2)=\zeta_1(y,z)+\zeta_2(x,y)=\zeta(x,y,z)(1+y^2)
	\quad\text{for $0<x,y,z< \sqrt{R_{e}}$.}
\end{align}
	
At first we examine the function $\zeta(x,y,z).$
It is not difficult to show that
 \begin{align*}
 	\zeta_2(x,y)=-\frac{10 x^3}{7 \left(x^2+1\right)}-\frac{x^3 y^2(10-7y)}{7 \left(y^2+1\right)}
 	-\frac{(10y+7)y^3}{7 \left(y^2+1\right)}
\frac{x^3}{1+y^2+y^2x^2}.
 \end{align*}
Since $x^3/(c_1+c_2x^2)$ ($c_1>0$ and $c_2\ge0$ are constants) is increasing with respect to $x$,
the function $\zeta_2$ and $\zeta$ are decreasing with respect to $x$ so that
\begin{align}\label{lemproof: appendix-DGS-condition2-x}
	\zeta(x,y,z)>\zeta(\sqrt{R_{e}},y,z)\quad\text{for $0<x,y,z< \sqrt{R_{e}}$}.
\end{align}
By following the elementary arguments in \eqref{lemproof: appendix-DGS-condition1-1}-\eqref{lemproof: appendix-DGS-condition1-3}, one can verify that $\frac{\partial}{\partial y}\zeta(\sqrt{R_{e}},y,z)<0$;
but the tediously long detials are omitted here. As depicted in Figure 1 (b),
$\zeta(\sqrt{R_{e}},y,z)$ has no any extreme points in $(0,R_3)^2$.
We consider the minimum value along the four boundaries:
\begin{itemize}
	\item [(i)] Along the side $y=0$, we have
	$\zeta(\sqrt{R_{e}},0,z)=2 -10 R_e\sqrt{R_e}/(7 +7R_e)\approx0.9579$.
	\item [(ii)] Along the side $y=\sqrt{R_e}$,
	\begin{align*}
	\zeta(\sqrt{R_{e}},\sqrt{R_{e}},z)\approx&\,
	\frac{-0.372002 z^5+0.0240772 z^4+0.638833 z^2+0.104165}{z^4+1.40198 z^2+0.401978}.		
	\end{align*}
    It has a unique maximum point at $z_3^*\approx0.448068$ and then
    	\begin{align*}
    	\zeta(\sqrt{R_{e}},\sqrt{R_{e}},z)
    	>&\,\min\left\{\zeta(\sqrt{R_{e}},\sqrt{R_{e}},0),
    	\zeta(\sqrt{R_{e}},\sqrt{R_{e}},z_3^*),\zeta(\sqrt{R_{e}},\sqrt{R_{e}},\sqrt{R_{e}})\right\}\\
    	=&\,\zeta(\sqrt{R_{e}},\sqrt{R_{e}},\sqrt{R_{e}})>\frac{1}{50}
    	\qquad\text{for $0< z< \sqrt{R_{e}}$.}		
    \end{align*}
	\item [(iii)] Along the side $z=0$, we have
	\begin{align*}
		\zeta(\sqrt{R_{e}},y,0)\approx&\,
		\frac{1.11457 y^5-0.676283 y^4-0.281384 y^3+1.105 y^2+0.385086}
		{\left(y^2+1\right)^2 \left(y^2+0.401978\right)}.		
	\end{align*}
It is decreasing such that
$\zeta(\sqrt{R_{e}},y,0)>\zeta(\sqrt{R_{e}},\sqrt{R_{e}},0)\approx0.259131$
for $0< y< \sqrt{R_{e}}$.
	\item [(iv)] Along the side $z=\sqrt{R_{e}}$, one has
	\begin{align*}
		\zeta(\sqrt{R_{e}},y,\sqrt{R_{e}})\approx&\,
		\frac{\left(y^2+1\right)^{-2}}{y^4+2.07416 y^2+0.672179}
		\Big(
		-0.790618 y^9-1.48448 y^7+1.32372 y^6\Big.\\
		&\,\Big.-0.825252 y^5+2.77809 y^4-1.06972 y^3+3.03679 y^2+0.643932\Big).
	\end{align*}
	It is also decreasing such that
	\begin{align*}
		\zeta(\sqrt{R_{e}},y,\sqrt{R_{e}})>\zeta(\sqrt{R_{e}},\sqrt{R_{e}},\sqrt{R_{e}})>\frac{1}{50}
		\quad\text{for $0<y< \sqrt{R_{e}}$.}		
	\end{align*}
\end{itemize}
It follows from \eqref{lemproof: appendix-DGS-condition2-x}
that $\zeta(x,y,z)>\frac{1}{50}$ for $0<x,y,z< \sqrt{R_{e}}$.
According to \eqref{lemproof: appendix-DGS-condition2-0},
we have $p(x^2,y^2,z^2)>\frac{1+y^2}{50}$ for $0<x,y,z< \sqrt{R_{e}}$.
This completes the proof.
\end{proof}

\end{document}